\documentclass[a4paper,11pt]{article}

\usepackage{color}

\usepackage[english]{babel}
\usepackage[latin1]{inputenc}
\usepackage{lmodern}
\usepackage{float}
\usepackage[T1]{fontenc}

\usepackage[a4paper,top=3cm,bottom=3cm,left=2cm,right=2cm,
bindingoffset=5mm]{geometry}

\usepackage{authblk}
\usepackage{amsmath,amssymb,amsthm}
\usepackage{graphicx}
\usepackage{multirow}

\newtheorem{remark}{Remark}

\newcommand\dd{\,{\rm{d}}}

\newcommand\M{{\bf M}}

\newcommand\Po{{\mathbb P}}
\newcommand\Pit{\widetilde P}
\newcommand\R{{\mathbb R}}

\newcommand\BB{{\mathcal B}}

\renewcommand\SS{{\cal  S}}
\newcommand\SSS{{\eta}}

\newcommand\Th{{\mathcal T}_h}
\newcommand\TT{{\mathcal T}}
\newcommand\TTT{\widetilde{\TT}}

\newcommand\gu{{\textsc{u}}}
\newcommand\gv{{\textsc{v}}}

\newcommand\MM{{\cal C}}

\newcommand\ff{\rm F}

\begin{document}

\centerline{\Large \bf Virtual Elements on polyhedra with a curved face }

\bigskip
\begin{center}
{\large{Franco Brezzi} \\
{\small \it IUSS, Piazza della Vittoria 15, 27100 Pavia,\\
 and IMATI del CNR, Via Ferrata 1, 27100 Pavia, (Italy)}}
\end{center}

\bigskip
\begin{center}
{\large{L. Donatella Marini} \\
{\small \it Dipartimento di Matematica, Universit\`a di Pavia,\\
 and IMATI del CNR, Via Ferrata 1, 27100 Pavia, (Italy)}}
\end{center}

\begin{abstract}
We revisit classical Virtual Element approximations on polygonal and polyhedral decompositions. We also recall the treatment proposed for dealing with decompositions into polygons with curved edges. In the second part of the paper we introduce a couple of new ideas for the construction of VEM-approximations on domains with curved boundary, both in two and three dimensions. The new approach looks promising, although sound numerical tests should be made to validate the efficiency of the method.
\end{abstract}


\section{Introduction}
The aim of this paper is twofold. On the one hand, it is an attempt to give an idea of the nature of the {\it Virtual Element Methods} (VEM) that were introduced around ten years ago  for the numerical solution of PDEs. This part would be essentially aimed at researchers that are {\it not}  specialists in Numerical Analysis and Scientific Computing.  There, on a polygonal domain $\Omega$ in 2d (polyhedral in 3d) we will recall the classical VEM approach. In the second part we will  first recall the approximation presented in \cite{Curvi2d} for dealing with 2d problems on a domain with a curved boundary. Then we will introduce new ideas which can simplify the previous approach, and, at the same time, allow the extension to 3-dimensional problems where the computational domain has a boundary (or part of it) described by a curved surface.

 Note that,  in many cases, an approximation of the boundary by a polyhedral surface would be too rough for the applicative point of view, unless the mesh is very fine (which, in three dimensions, could be  outrageously expensive). 
  Moreover, the approximation of the curved boundary with  straight lines in 2d (flat polygons in 3d) might produce a loss of convergence. This occurs with Finite Element approximations, and is known as Bab\v uska paradox (see  \cite{Babuska71}). Since on triangular decompositions Virtual Elements of degree $k=1$ coincide with piecewise linear Finite Elements, it is highly possible that a similar situation might happen with VEM, even in more general cases. Thus the interest of devising VEM-approximations able to deal with curved boundaries (see, e.g., \cite{VEM-curv-Bertol}, \cite{rozzi},  \cite{Ovall}, \cite{Curvi2d}).
 
The possibility of dealing with curved boundary is one of the main reasons that ensured the enormous success of the so-called  IsoGeometric Analysis (IGA; see e.g. \cite{TJR-Iso}), where the solution is approximated by splines of the same type used  (in almost 100\% of industrial design) to describe the boundary of the domain. 
 
 Here, however, contrary to IGA, we do not need to correlate the two ways of {\it describing the boundary} and {\it discretizing the problem}. In some sense, the type of {\it description of  the boundary} and the type of {\it discretization of the problem} are, in essence, 
{\it totally independent one of the other} and the two choices do not need to be correlated. This, obviously, does not mean that you can choose the two discretizations independently, but that you can choose the {\it type of discretization} (the decomposition and the degree to be used in each polyhedron), essentially, without considering
the {\it type of description} of the boundary: by splines, trigonometric functions, or just (brutally) by a given analytic expression.

We will use the Poisson problem as super-simple example in order to explain the construction of the numerical method (passing from the PDE to a system of linear equations that can then be solved on the computer).

The problem of constructing three-dimensional Virtual Element spaces on  general polyhedra
 having one or  more curved faces  is still open (at the fully general level). Here we propose several possible ways of tackling the problem.

We point out that the ``discretizations of order $k$" described here will satisfy the so-called {\it patch-test of order $k$} (very popular among Engineers) that, very roughly, says that: {\it whenever the exact solution is a polynomial of degree $k$, then the solution of the discretized problem will coincide with the exact solution}.
Although the patch-test is only a necessary, but not a sufficient condition for convergence, it makes you feel more confortable to have it satisfied.

 From a certain point of view, insisting on the correct treatment of curved boundaries,  the approach  presented here could constitute a viable alternative to the classical Iso-Geometric Analysis (already mentioned), allowing a wider choice of the domain decomposition, payed with a lower regularity of the subspace (basically $H^1$ instead of the $H^p$  regularity of IGA). Note that, compared with Finite Elements, the methods proposed here
would indeed allow a much easier use of more regular subspaces (say, like $H^2$ or $H^3$) although using more degrees of freedom than IGA or, more generally, than splines. 

An outline of the paper is as follows.

After a short reminder of classical notation, we will start recalling classical VEM spaces on polygons/polyhedra, and their use for the approximate solution of PDE problems. 

 This will include both the definition of the VEM spaces, and a hint on the classical ways of using them to solve PDE problems on the computer.  We anticipate that, within each element of the decomposition, the VEM spaces, in general, contain (together with all polynomials of a certain degree) also  functions that are {\it not} polynomials, and whose point values are not known. Typically, in each element  they will be solutions of suitable boundary value problems whose boundary data and right-hand side are identified by a certain number of parameters (the {\it degrees of freedom}). The solution of such element-wise boundary value problems is clearly out of the question, and each element of the Virtual Element space will be known only through its degrees of freedom. This implies that its manipulation will not be as easy as it happens when using just polynomials (or splines), and we will briefly indicate how all this can be dealt with on the computer.

Then we will discuss the use of  {\it polygons with a curved edge}, recalling what was done in \cite{Curvi2d} and proposing some new interesting alternatives. At this point we will be ready to discuss the treatment of polyhedra with a curved face. For simplicity we will deal only with the case in  which each element has just {\it one}  curved face,  that is  part of the boundary of the computational domain, but the extension to more general cases should not be too difficult.

We believe that the approaches that we suggest might open new perspectives on the treatment of curved edges/faces, and lead  to new more suitable methods. However, sound numerical tests will have to be performed to validate the efficiency and, most important, the accuracy of the methods, but this goes beyond the  scopes of the present paper.

\subsection{Notation}

Throughout the paper we will use the common notation for functional spaces $L^2(D)$ and $H^s(D)$ on a domain $D$, with scalar product and norm $(\cdot,\cdot)_{0}, \|\cdot\|_0$ ( and $(\cdot,\cdot)_{s}, \|\cdot\|_s$, respectively); $H^{1/2}(\Gamma)$ will denote the space of traces of functions of $H^1(D)$ on the boundary $\Gamma$  (see  \cite{Lions:Magenes:1968}).

When convenient, we will use the notation $\Delta_2$ and $\Delta_3$ to indicate the two-dimensional or the three-dimensional Laplacian, respectively.

For $k\ge 0$, $\Po_{k}(D)$ will denote the space of polynomials of degree $\le k$ on a domain $D$. As common,  $\Po_{-1}=\{0\}$. 

Throughout the paper, $C$ will 
represent, as usual, a fixed positive constant, not necessarily the same from one occurrence to another.

We also point out that, given a domain $\Omega$ and a subset $P$ of it, by saying that "$P$ is {\it internal  to $\Omega$}"  we mean that the closure of $P$ is a subset  of  $\mathop \Omega\limits^ \circ$.

Given a finite dimensional linear space $V$ and an integer number $S$,  a linear mapping 
\begin{equation}\label{defGen}
{\mathcal G}\in\mathcal L(\R^S,V)
\end{equation}
  is said to be  {\it a set of generators for} $V$ if it is {\it surjective}. If moreover ${\mathcal G}$ is also {\it injective}, then it is one-to-one,
and its inverse $ {\mathcal D}: V\rightarrow\R^S$ could be represented as a set of $S$  linear operators 
 $(\delta_1,...,\delta_S)$, each in $\mathcal L(V,\R)$, that are then called {\it degrees of freedom}:
\begin{equation}\label{defD}
{\mathcal D}:\, v\rightarrow (\delta_1(v), \delta_2(v),...,\delta_S(v))  (= {\mbox{\it  the degrees of freedom of }  v}).
\end{equation} 
Note that, even when ${\mathcal G}$ is not injective (but only surjective) we could always construct a {\it right inverse}
${\mathcal D}: V\rightarrow\R^S$ such that
\begin{equation*}
{\mathcal G}{\mathcal D}v=v\quad \forall v\in V
\end{equation*}
that however, this time, will not be {\it unique}. 

In what follows, a set of generators for a space $V$ will be typically denoted by 
${\bf G}_V$.


\subsection{Projectors in VEM spaces}\label{projectors}
As we shall see, most of the Virtual Element spaces that we are going to construct will be made of functions that are solutions of local (i.e., element by element) boundary value problems for PDEs. Typically, in two dimensions, we will know explicitly only their values
at the interelement boundaries, and possibly some of their moments inside each element. The computation of other quantities, as for instance their pointwise value inside the elements, will be totally out of reach. 
Consequently, we will often (actually: almost always) use in their place some suitable {\it projection} on polynomial spaces. There too, from the information that we have, some projections will be computable, and others will not.
 For the moment we just briefly anticipate the classical structure of VEM spaces of order $k$ on a polygonal element $P$, in order to see some projectors that are actually computable out of the degrees of freedom.

The {\it general structure} of a VEM space with order of accuracy $k\ge 1$ on a polygon $P$ is:
\begin{equation}\label{typical}
V_k(P):=\{v\in C^0(\overline{P}) \mbox{ such that } v_{|e}\in \Po_k(e) ~\forall \mbox{ edge } e,
\mbox{ and }\Delta v\in \Po_{k-2}(P)\}.
\end{equation}
Typically the degrees of freedom for the spaces in \eqref{typical} will provide directly
\begin{itemize}
\item
the values on each edge of $P$, and 
\item (for $ k\ge 2$)  the moments up to the order $k-2$ in $P$.
\end{itemize} 

\begin{remark} An immediate, important, generalization of the  spaces \eqref{typical} can be obtained separating the degree of the functions on each edge and the degree of the Laplacian. Hence, say, for $k\ge 1$ and $k_L \ge -1$ we can consder the spaces
\begin{equation}\label{typical2}
V_{k,k_L}(P):=\{v\in C^0(\overline{P}) \mbox{ such that } v_{|e}\in \Po_k(e) ~\forall \mbox{ edge } e,
~\Delta v\in \Po_{k_L}(P)\}.
\end{equation}
where the degrees of freedom for the spaces in \eqref{typical2} should provide directly
\begin{itemize}
\item
the values on each edge of $P$, and 
\item  the moments up to the order $k_L$ in $P$.
\end{itemize} 
Clearly the spaces \eqref{typical} would be obtained for $k_L=k-2$.
\qed
\end{remark}

Here below we will see some computable projections that we are going to use in the sequel.

\medskip
\noindent {\bf *  The $H^1_0(P)$-projection}.  Given a polygonal element $ P$, for each $v\in H^1(P)$ we define its  projection  $\Pi^{\nabla}_kv$ onto the space $\Po_k$
as the solution, in $\Po_k(P)$, of
\begin{equation}\label{defpinabla}
\Pi^{\nabla}_k v \in \Po_k(P), {\mbox{ and}} \int_P \nabla(\Pi^{\nabla}_k v)\cdot\nabla q_k\,\dd P=\int_P \nabla v \cdot\nabla q_k\,\dd P\quad \forall q_k\in \Po_k.
\end{equation}
Actually, \eqref{defpinabla} identifies $\Pi^{\nabla}_kv$ only up to a constant, that can easily be fixed, for instance,  with the additional requirement that
\begin{equation}\label{fixc}
\int_{\partial P} \Pi^\nabla_k v\,\dd s=\int_{\partial P} v\,\dd s.
\end{equation}
The left-hand side of \eqref{defpinabla} is a product of polynomials, and is obviously computable. Integrating the right-hand side by part we have 
\begin{equation}\label{byparts}
\int_P \nabla v \cdot\nabla q_k\,\dd P=-\int_P  v\,\Delta q_k\,\dd P + \int_{\partial P} v \frac{\partial q_k}{\partial n}\,\dd s
\end{equation}
and we will be able to perform our computation as far as we know explicitly
\begin{equation*}
\left\{
\begin{aligned}
&\mbox{the moments of $v$ up to the order $k-2$ on $P$},\\
&\mbox{the moments of $v$ up to the order $k-1$ on each edge of $P$}.\\
\end{aligned}
 \right.
\end{equation*}
\noindent As we anticipated, normally on each edge we  actually know directly the whole (polynomial) values of the functions of our discretized spaces, and the integrals on the edges can be computed exactly.
Moreover, the first term on the right-hand side of \eqref{byparts} can be easily computed out of the degrees of freedom.

 Now let us see some other  important projectors.

\medskip
\noindent {\bf * The $L^2$ projection of the gradient}. With an almost identical argument we can compute the $L^2$-projection $\Pi_{k-1}^0(\nabla v)$ on $[\Po_{k-1}]^2$, defined by
\begin{equation}\label{projgrad}
\int_P (\Pi_{k-1}^0\nabla  v) \cdot {\bf q}\,\dd P=\!\!\int_P \nabla v \cdot{\bf q}\,\dd P 
\,\qquad\forall {\bf q}\in [\Po_{k-1}]^2.
\end{equation}
Indeed, integrating by parts the right-hand side we have
\begin{equation*}
\int_P \nabla v \cdot{\bf q}\,\dd P=\int_{\partial P} v\, {\bf q }\cdot {\bf n}\dd s-\int_P v \,{\rm div}{\bf q}\dd P ,
\end{equation*}
and both terms are immediately computed as in \eqref{byparts}.
Note that, referring to  \eqref{typical2}, we can actually compute the $L^2$ projection of $\nabla v$ on $(\Po_s)^2$ whenever
$v\in V_{k,k_L}$ with $k_L\ge s-1$.

\medskip
\noindent {\bf * The {\it dofi-dofi} projector}.
It is also important  to note that there are other projectors from $V_k(P)$ to $\Po_k$ that  are computable out of the degrees of freedom ${\mathcal D}v$. The simplest one, that we call $\Pi^{\mathcal D}$, would be defined for each $v\in V_k(P)$ as the (unique) solution in $\Po_k$ of
\begin{equation}\label{proD}
\Big({\mathcal D} (\Pi^{\mathcal D} v)-{\mathcal D}v,{\mathcal D}(q_k)\Big)_{\R^S}=0\quad\forall q_k\in\Po_k,
\end{equation}
where obviously $(\cdot\,  ,\cdot)_{\R^S}$ is the usual Euclidean scalar product in $\R^S$, and $S$ is the number of degrees of freedom of $V_k(P)$.   Note that this can be done even if in \eqref{proD} the mapping ${\mathcal D}$ is not representing the degrees of freedom, but just
any other {\it identifier} of the type \eqref{defD}, provided that ${\mathcal D}$ is injective from $\Po_k$ to $\R^S$. Indeed, as far as $S$ is equal to the number of degrees of freedom,
\eqref{proD} is just a pompous way to say that ${\mathcal D} (\Pi^{\mathcal D} v)={\mathcal D}v$.
\qed

\medskip
\noindent{\bf * The $L^2$ projection.}  
In general, the $L^2$- projection of an element $v$ on $\Po_s(P)$ can be computed only when we know the moments of $v$ {of order up to $s$}.
If we know the integrals $\int_P v m_j$ for $ j=1,...S$, where the $m_j$ are a basis for $\Po_s$, and $S$ is the dimension of $\Po_s$,  then we can orthonormalize the $m_j$. Denoting by $\widetilde{m_j}$ the orthonormal basis, i.e., such that
\begin{equation*}
\qquad \int_P \widetilde{m_j}\widetilde{m_i}\dd P=\delta_{i,j} \qquad i,j=1,...,S,
\end{equation*}
the $L^2$ projection $\Pi^{0,P}_s v$ will be
\begin{equation*}
{\Pi^{0,P}_s v}:=\sum_{j=1}^{S} \widetilde{m}_j{\int_P  v \widetilde{m}_j\dd P}. 
\end{equation*}

\begin{figure}[!htb]
\begin{center}
\includegraphics[height=0.22\textwidth]{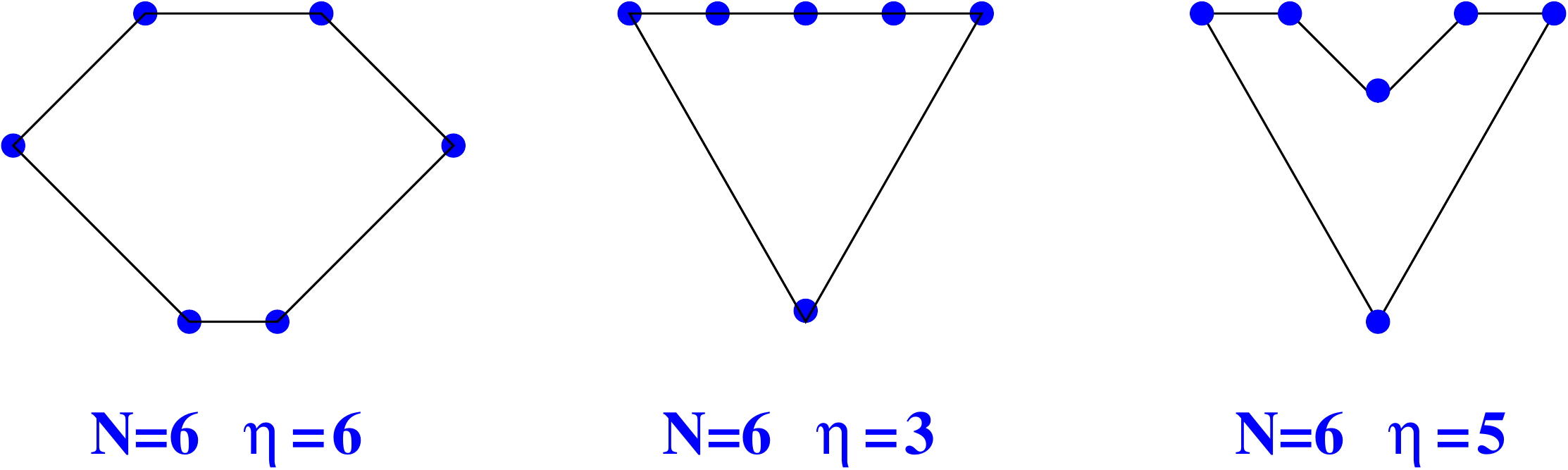} 
\end{center}
\caption{The value of $\eta$ for various polygons with $N$ edges}
\label{Vari eta}
\end{figure}

\medskip
\noindent {\bf * The Serendipity-like projectors}.
 A particularly relevant  class of possible alternative projectors is given by the {\it Serendipity-like} projectors. The basic idea (or, so to speak, {\it the seed} of it)
is the following:  For an element $P$ we denote by $\SSS(P)$ {\it the minimum number of straight lines necessary to cover the whole boundary $\partial P$}. Note that in our setting two edges (consecutive or not) may lie on the same straight line (see an example in Fig. \ref{Vari eta}). Then we observe that for every $k <\SSS(P) $ a polynomial of degree $k$ will be completely identified by its value on the boundary $\partial P$.
Then, always for $k<\eta(P)$   we can define a projector $\Pi^\SSS_k$ from $V_k(P)$ to $\Po_k(P)$ 
defined by
\begin{equation}\label{projbordo1}
\int_{\partial P}(v-\Pi^\SSS_k v)q_k\dd s=0\quad \forall q_k\in \Po_k(P)
\end{equation}
which clearly has a unique solution. For $k\ge \SSS(P)$ we will have to add some {\it internal} information. Typically, we may choose an integer $r$ with $k-\SSS(P)\le r\le k-2$, add to \eqref{projbordo1} the condition
\begin{equation}\label{projdentro1}
\int_{ P}(v-\Pi^\SSS_k v)q_r\dd s=0\quad \forall q_r\in \Po_r(P),
\end{equation}
and consider the {\it pair} of equations \eqref{projbordo1}-\eqref{projdentro1}, possibly to be solved in the {\it least squares} sense,
mimicking what is done in the construction of Serendipity VEM spaces (see \cite{SERE-nod}). 

An important remark is in order, concerning the $L^2-$ projection of the gradient, in the context of the {\it stabilization procedures} that will be discussed in the next pages.

\begin{remark} In several cases, it would be convenient (and we daresay {\it very convenient}) to construct polynomial projections   $\Pi (\nabla v)$ of the gradient that satisfy
\begin{equation}\label{iniettiva}
\Pi\nabla  v=0 \quad =>\quad \nabla v=0 
\end{equation}
 (or, in oher words, that are {\it injective}). Clearly, if the polygonal element that we are considering has {\it many
	edges} (and consequently, even for a low $k$, {\it many degrees of freedom}), condition \eqref{iniettiva} would require
that we project onto  a {\it big} polynomial space, and will make the whole idea {\it too expensive}. However, as we shall see in a while, the advantages of  \eqref{iniettiva} are relevant in many cases, and the decision on {\it whether to enforce it or not} could be delicate in several circumstances.
\qed
\end{remark}

\subsection{The simplest model problem}\label{disproret}
Let $\Omega\subset \R^{d}, d=2,3,$ be a domain with boundary $\Gamma$. We assume that $\Gamma$ is the union of two (for simplicity, connected) parts, denoted by $\Gamma_D$ and $\Gamma_N$. As usual we will assume that $\mathring\Gamma_D\cap
\mathring\Gamma_N=\emptyset$ and $\overline\Gamma_D\cup\overline\Gamma_N\equiv\Gamma$.  Let $f, g_D, g_N$ be given functions with $f\in L^2(\Omega)$, $g_D\in H^{1/2}(\Gamma_D$), and $g_N$, say, in $ L^2(\Gamma_N)$. We consider the simple model problem
\begin{equation}\label{Pbmod}
\left\{
\begin{aligned}
\mbox{find } &u\in H^1(\Omega) \mbox{ such that}\\
-\Delta u &= f \quad \,\mbox{  in } \Omega,\\
u &= g_D \quad \mbox{on } \Gamma_D,\\
\frac{\partial u}{\partial n} &= g_N \quad \mbox{on }\Gamma_N.
\end{aligned}
 \right.
\end{equation}
Setting, for $\varphi\in H^{1/2}(\Gamma_D)$,
\begin{equation*}
H^1_{\varphi,\Gamma_D}(\Omega):=\{v\in H^1(\Omega)\mbox{ such that } v=\varphi \mbox{ on }\Gamma_D\},
\end{equation*}
the variational formulation of \eqref{Pbmod} can be written as
\begin{equation}\label{Pbmod2}
\left\{
\begin{aligned}
&\mbox{ find } u\in H^1_{g_D,\Gamma_D}(\Omega) \mbox{ such that:}\\
&a(u,v)=(f,v)+ <g_N,v>_{\Gamma_N}\quad\forall v\in H^1_{0,\Gamma_D}(\Omega)
\end{aligned}
\right.
\end{equation}
where $a(u,v)$ is the bilinear form defined by
\begin{equation}\label{defauv}
a^P(u,v)=\int_{P}\nabla u\cdot\nabla v\;\dd P,\qquad  a(u,v)=\sum_P a^P(u,v),
\end{equation}
and
\begin{equation}
(f,v)=\int_{\Omega} f\,v \dd \Omega, \quad <g_N,v>_{\Gamma_N}=\int_{\Gamma_N} g_N\, v \dd \Gamma,
\end{equation}
and the boundary integral could be replaced by  a suitable duality for a less regular  $g_N$.

For the sake of simplicity, here we will discuss, separately, only the two cases: {\it Full Dirichlet} (when $\Gamma_D\equiv\Gamma$; {actually, very simple)} and {\it Full Neumann} (when $\Gamma_N\equiv\Gamma$). 

In both cases, we will assume that we are given a sequence of decompositions $\{\Th\}_h$ of $\Omega$  in {polytopes} $P$ with diameter $h_P$, and we indicate with $|h|$ the maximum of $h_P$ for $P\in \Th$. We will also make the usual assumptions that each polytope (of each $\Th$) is star-shaped with respect to a ball of radius $\rho_P\ge C h_P$, and each edge/face of a polytope $P$ has length $\ge C h_P$.

\section{Virtual Elements for 2D polygons}\label{2dflat}

\subsection{The local VEM spaces}\label{local-space}
For the sake of simplicity, we  begin by recalling the original {\it plain vanilla} VEM spaces on two-dimensional polygons with straight edges, as presented in \cite{volley} and already anticipated in \eqref{typical}. 
Let $P$ be a polygon, and let $k$ be an integer $\ge 1$.

We define first
\begin{equation}\label{defBkf}
\BB_k(\partial P):=
\{v\in C^0(\partial P) \mbox{ s.t. } v_{|e} \in \Po_k(e)\, \forall \mbox{ edge } e \in \partial P\},
\end{equation}
and then
\begin{equation}\label{Vkf}
V_k(P)=\{v\in H^1(P) \mbox{ s.t. } v_{|\partial P} \in \BB_k(\partial P) \mbox{ and }\Delta v \in \Po_{k-2}(P)\}.
\end{equation}

It is clear that all polynomials of degree $\le k$ belong to $V_k(P)$, and that (as anticipated in the previous section)   an element $v$  of $V_k(P)$ is uniquely determined  by:
\begin{align}
&\mbox{the values of }v \mbox{ at the vertices of  }P,\label{dofs-face-vert}\\
&\label{dof-face-edg}
\mbox{(for }k\ge 2)\quad\int_e v\, p_{k-2}\dd e\quad\forall p_{k-2}\in \Po_{k-2}(e)\quad\forall\mbox{ edge }e \subset  \partial P,\\
&\label{dof-face-ins}
\mbox{(for }k\ge 2)\quad\int_P v\, p_{k-2}\dd P\quad\forall p_{k-2}\in \Po_{k-2}(P).
\end{align}
For reasons that will be clear soon, the degrees of freedom \eqref{dofs-face-vert}--\eqref{dof-face-ins} should be suitably scaled so that they all scale in the same way. We do not enter the details here, for which we refer for instance to \cite{hitchhikers}.

\subsection{The  two-dimensional ``full Dirichlet" case}\label{disproretD}
Let $\Omega\subset \R^{2}$ be a polygonal domain with boundary $\Gamma$.  Let moreover $g$ and $f$  be smooth-enough functions defined on $\Gamma$ and in $\Omega$, respectively.
Our model problem \eqref{Pbmod2} becomes
\begin{equation}\label{PbmodD}
\left\{
\begin{aligned}
\mbox{find } &u\in H^1(\Omega) \mbox{ such that}\\
-\Delta u &= f \quad \,\mbox{  in } \Omega,\\
u &= g \quad \mbox{on } \Gamma.
\end{aligned}
 \right.
\end{equation}
Since now $\Gamma_D\equiv\Gamma$,  we can simply write  $H^1_g(\Omega)$ instead of $H^1_{g,\Gamma_D}(\Omega)$, and \eqref{Pbmod2} becomes
\begin{equation}\label{Pbmod-Dir}
\left\{
\begin{aligned}
&\mbox{ find } u\in H^1_{g}(\Omega) \mbox{ such that:}\\
&a(u,v)=(f,v)\quad\forall v\in H^1_{0}(\Omega) .
\end{aligned}
\right.
\end{equation}
For  a generic $\varphi$ defined on $\Gamma$ (here either $\varphi=g$ or $\varphi=0$), and for an element $P\in\Th$, we define the local spaces as
\begin{equation}\label{VkfP}
\begin{aligned}
V_k^\varphi(P)&=\{v\!\in\! H^1(P)\mbox{ s.t.}\!: v_{|e}\in\Po_k(e)\;   \forall  e\notin \Gamma,\;
 v_{|e}=\varphi_{|e} \; \forall e \in { \Gamma}, \\
 &\hskip2cm \mbox{and   } \Delta v \in \Po_{k-2}(P)\},
 \end{aligned}
\end{equation}
with their global counterpart
\begin{equation}\label{VkfOm}
V_h^\varphi(\Omega)=\{v\in H^1(\Omega) \mbox{ s.t. }   v_{|P}\in V_k^\varphi({P}) \;\forall P\in\Th \}.
\end{equation}
The degrees of freedom in $V_h^\varphi(\Omega)$, on top of the value of $\varphi$ on $\Gamma$, will obviously be:
\begin{itemize}
\item The values at the internal vertices of the decomposition,
\item (for $k\ge 2$) The moments on each internal edge up to the degree $k-2$,
\item (for $k\ge 2$) The moments inside each element up to the degree $k-2$.
\end{itemize}
It is immediate to check that (as seen in the subsection \ref{projectors}) we can follow \eqref{defpinabla} and \eqref{fixc} and define, in each element $P$, the projector
$\Pi^{\nabla}_k$, that will be computable using the above degrees of freedom.

Next, always on each element $ P$, we consider the restriction  $a^P(u,v)$ of the bilinear form $a(u,v)$  to $P$, and then on each $V_k^\varphi(P)$ (as defined in \eqref{VkfP}) we set
\begin{equation}\label{defahpd}
a^P_h(u,v):=a^P(\Pi_k^\nabla u, \Pi_k^\nabla v)+\SS^P(u-\Pi_k^\nabla u,v-\Pi_k^\nabla v) .
\end{equation}
The first term in the right-hand side of \eqref{defahpd} is the consistency term which is in general singular, and the second term is needed for stability. We can take as  $\SS^P(u,v)$  any bilinear form that scales as $a^P(u,v)$ and is positive on the kernel of $\Pi^\nabla_k$. In particular we would require that there exist two positive constants $\alpha_*$  and  $\alpha^*$ such that
\begin{equation*}
\alpha_* a^P(v_h,v_h)\le a^P_h(v_h,v_h)\le\alpha^*a^P(v_h,v_h) \qquad \forall v_h \in V_k^{\varphi}(P).
\end{equation*}

 \noindent The simplest example would be {(setting $S=$ number of dofs in $P$)}: 
\begin{equation*}
\SS^P(u,v):=\sum_{i=1,S}
do\!f_i(u)do\!f_i(v),
\end{equation*}
where $dof_i$ is the $i^{th}$ degree of freedom properly scaled in such a way that $\SS^P(u,v)$
scales as $a^P(u,v)$. This is the reason why the degrees of freedom need to be properly scaled, as briefly anticipated at the end of Sect. \ref{local-space}.
But other choices can be convenient, as 
\begin{equation*}
\SS^P(u,v):=h_P^{-1}\int_{\partial P} u\,v\dd s ,
\end{equation*}
{(where $h_P$ is still the diameter of the element $P$),} or
\begin{equation*}
\SS^P(u,v):=h_P\int_{\partial P} u_t\,v_t\dd s ,
\end{equation*}
where $u_t$ and $v_t$ are the {\it tangential derivatives} of $u$ and $v$, respectively.
Needless to say,  the above expressions 
for $\SS^P(u,v)$ can be multiplied by a fixed constant (say, 5, or 1/5, or whatever fixed constant you might like, as far as it is independent of $h$).

\begin{remark} Instead of using the $\Pi^{\nabla}_k$ operator in the definition of the consistency part of the discrete bilinear form \eqref{defahpd}, we might use the $L^2$-projection of the gradient, and define
\begin{equation}\label{variable}
a^P_h(u,v):=\int_P \Pi^0_{k-1} (\nabla u)\cdot \Pi^0_{k-1} (\nabla v)+\SS^P( u-\Pi^{\nabla}_k u, v - \Pi^{\nabla}_k v) .
\end{equation}
Actually, this choice is preferable for more general problems, in particular in presence of variable coefficients, since 
there the $\Pi^ \nabla_k$ operator might produce a loss of order of convergence for high values of $k$ (see \cite{variable-primal}). \qed
\end{remark}

\begin{remark} For the case that we are considering here, most choices of the stabilizing form $\SS^P$ will give quite satisfactory results. However, for more complex problems the choice of the stabilizing term (or, for a given term, the choice of 
the   coefficient to be used in front of it) might become quite delicate: a poor choice would not jeopardize the convergence of the method (in the limit for $h$ going to $0$), but would give poor results for the {\it reasonably affordable decompositions}. \qed
\end{remark}
\begin{remark}  In some sense, we might  say that, in several applications (just to take an example: for  the extension of the present approach to fourth order problems), the choice of the right stabilizing term (with the right coefficient) can become {\it the most delicate point} (see \cite{Brezzi:Marini:plates}, \cite{Chinosi:Marini:plates}). Clearly a suitable amount of {\it good new ideas}  are required on this subject.  Recently, attempts of getting rid of the stabilizing form have been performed (separately, both by {\it Berrone and co-workers} and  by {\it Perego  and co-workers}: personal communications). The idea is to use the projection of the gradient, as in \eqref{variable}, onto polynomials of higher degree, for example onto $[\Po_k]^2$ instead of $[\Po_{k-1}]^2$ with the hope to produce  a nonsingular consistency part. So far a wide theoretical study is still to be made, but the preliminary experimental results are encouraging and show that this strategy is worth analizing, in particular for polygons having a number of edges {\it not too big}.
\qed
\end{remark}

Collecting  the definitions \eqref{defahpd}  we set
\begin{equation*}
a_h(u,v):=\sum_{P\in \Th} a^P_h(u,v)\qquad \forall u,v\in V^{\varphi}_h(\Omega).
\end{equation*} 
For a detailed treatment of the right-hand side we refer to \cite{volley}. Here we just recall the definition of $f_h$ on each polygon. With $\Pi^0_k:=L^2-$projection operator onto $\Po_k$, we set 
\begin{equation}\label{terminenoto}
{f_h}_{|P}=
\begin{cases}
\Pi^0_0 f \quad &\mbox{for } k=1,\\
\Pi^0_{k-2}f \quad &\mbox{for } k\ge 2.
\end{cases}
\end{equation}

\noindent The discretised problem reads now:
\begin{equation}\label{Pbmodh2d}
\left\{
\begin{aligned}
&\mbox{Find $u_h\in V_h^{g}(\Omega)$ such that}:\\
& a_h(u_h,v_h)=(f_h,v_h)
\quad\forall v_h \in V_h^0(\Omega).
\end{aligned}
\right.
\end{equation}

\begin{remark}
We point out that, if the boundary value $g$ is {\bf not} the trace of a polynomial, we will loose the property that {\it every local space
$V_k^g(P)$ contains all the polynomials of degree $\le k$ on $P$}. If, for some reason, we are  interested in such property, then
we should take first a piecewise polynomial approximation $g_h$ of $g$, and then  use $V^{g_h}_k(P)$ and $V^{g_h}_h(\Omega)$
in place of  $V^{g}_k(P)$ and $V^{g}_h(\Omega)$, respectively.
\qed
\end{remark}

\subsection{A tiny historical note}\label{mix}

As explicitly claimed from the very beginning (see e.g. \cite{volley}) the Virtual Element Methods are a direct follow-up (and initially almost a simple re-formulation) of Mimetic Finite Differences (MFD, see e.g. \cite{BLMbook} and the references therein). In MFD the unknowns are just the degrees of freedom (as, indeed,  it always happens in the final computer code) and the formulation is done directly in terms of them  (typically, nodal values or averages); in other words, with MFD we {\it do not} have a discrete functional space playing the role that here (and generally in VEMs) is played by $V_h$. There too, as in VEM, the contribution of each element to the final stiffness matrix, for instance for a problem as \eqref{Pbmod}, is the sum of a consistency part and  of a stabilizing part. The first part is obtained considering first the {\it polynomial  subspace} of dofs,  meaning the subspace whose elements are obtained, each,  as the {\it values at nodes} and  {\it the averages}  of a single polynomial of the prescribed degree $k$. So, for $k=1$ we will have the subset (of dimension 3) generated by the (say) nodal values of $1,x,y$. The consistency part of the MFD stiffness matrix will then act only on this subspace, and there it will reproduce the corresponding  expected behaviour of differential operators and/or integrals that appear in the continuous bilinear form (as \eqref{defauv}).
The stabilization part, instead, will vanish on the above subspace, and will be symmetric-positive-definite on its complement (more precisely, in the complement of its kernel). In other words, we have the perfect analogue of \eqref{defahpd}(that does actually {\it coincide} with the one used in MFD for a judicious choice of the respective stabilizing parts).

The two methods, MFD and VEM, are just part of a pack of different approaches to the treatment of PDEs using decompositions of the computational domain into elements whose shape is more general than just simplexes and boxes.
Among the most ancient ones (together with MFD) another building block is surely made by Discontinuous Galerkin Methods
(see, e.g. \cite{ABCM}, or the most recent  \cite{Dole-Feista},\cite{DiPietro-Ern}). There the discrete spaces are just made of polynomials, and suitable tricks are needed in order to enforce some kind of {\it weak continuity} from one element to another. Many variants are related with them, from Hybridizable Discontinuous Galerkin methods (where, roughly speaking, a weak continuity is enforced via Lagrange Multipliers),
or Hybrid-High-Order and Weak Galerkin methods, where one works with {\it two piecewise polynomial spaces}: one at the boundary and one inside. See e.g. \cite{cockburn-IMU}, \cite{DiPietro-Ern}, \cite{Mu-Wang-Ye} and the references therein. See also 
\cite{Gong-Wu-Xu} for similar ideas.  
Needless to say, if necessary one could always use one type of discretisation in some elements and another one in other elements, as, in the end, the computer will always solve ``just a problem in $\R^n$" having as unknowns the degrees of freedom.

\subsection{The two-dimensional{ ``full Neumann"} case}
We consider now the {\it full Neumann} case, i.e., $\Gamma_N\equiv \Gamma$.  We recall that, in this case, $g_N$ and $f$ must satisfy the compatibility condition 
\begin{equation*}
 \int_{\Gamma}g_N\dd s=\int_{\Omega}f\dd x
\end{equation*}
due to the Gauss divergence theorem. We also recall that the solution $u$ will be determined only up to an additive constant, that in the computer code can be fixed, as common, just by fixing its value to be equal to zero at some point (usually, a vertex of $\Omega$). From now on, when discussing a {\it full Neumann problem}, we will implicitly assume that we have chosen, once and for all, a vertex of the discretization, and prescribed that the solution of the continuous problem $u$, and all the elements of the discrete subspaces, vanish there.

The discrete problem differs from the {\it full Dirichlet} case only in the right-hand side, upon performing a slight change in the definition of the local and global discrete spaces \eqref{VkfP} and \eqref{VkfOm}. By defining
\begin{equation*}
V_k(P)=\{v\!\in\! H^1(P)\mbox{ s.t.}\!: v_{|e}\in\Po_k(e)\;   \forall  \mbox{ edge }e,  \;\mbox{and   } \Delta v \in \Po_{k-2}(P)\},
\end{equation*}
\begin{equation*}
V_h(\Omega)=\{v\in H^1(\Omega) \mbox{ s.t. }   v_{|P}\in V_k({P}) \;\forall P\in\Th \},
\end{equation*}
the discrete problem reads:
\begin{equation}\label{Pbmodh2d-Neumann}
\left\{
\begin{aligned}
&\mbox{Find $u_h\in V_h(\Omega)$ such that}:\\
& a_h(u_h,v_h)=(f_h,v_h)+<g_N,v_h>_{\Gamma}
\quad\forall v_h \in V_h(\Omega).
\end{aligned}
\right.
\end{equation}
In \eqref{Pbmodh2d-Neumann}  $a_h(u_h,v_h)$ and $(f_h,v_h)$ are the same as in the previous subsection, and  $<g_N,v_h>_{\Gamma}$ is a boundary integral that can be computed since the functions $v_h$ are polynomials on each edge, completely known.

\section{ Recalling 3D VEM for classical polyhedra}\label{3dpiatti}

The extension of the above construction to the three dimensional case is, in some sense, both immediate and tricky. Indeed, the first attempt that comes to mind, given a polyhedron $P$ and an integer $k$, is to extend slavishly what we did in the two-dimensional case. Let us first examine the {\it full Dirichlet case}.

To begin with, we define 
\begin{equation}\label{PbmodD3d}
\begin{aligned}
V_{h}(\Omega)&:=\{v\in H^1(\Omega)\cap C^0(\overline{\Omega}) \mbox{ such that: }
v_{|e}\in \Po_k(e)\;\forall \mbox{ internal edge } e,\\
&\Delta_2 v_{|{\ff} }\in \Po_{k-2}(\ff)\,\forall \mbox{ internal face } \ff, ~
\Delta_3 v\in \Po_{k-2}(P)\, \forall \mbox{ polyhedron } P\}.
\end{aligned}
\end{equation}

We remark that $V_{h}(\Omega)$, as defined in \eqref{PbmodD3d}, is infinite dimensional (we are not making requirements on the values of $v$ on  $\partial \Omega$).
Then, given a (smooth enough) function $\varphi$ defined on $\partial\Omega$ we can mimic \eqref{VkfP}-\eqref{VkfOm} and restrict our space  to
\begin{equation}\label{3dVhkphi}
V^{\varphi}_{h}(\Omega):=V_{h}(\Omega)\cap H^1_{\varphi}(\Omega)
\end{equation}
that, now, will be finite-dimensional. In particular, for a fixed given $\varphi$, the degrees of freedom in \eqref{3dVhkphi} will be:
\begin{equation}\label{dof3d}
\left\{
\begin{aligned}
&\mbox{The values at the internal vertices of the decomposition},\\
&\mbox{(for $k\ge 2$) The moments on each internal edge up to the degree $k-2$},\\
&\mbox{(for $k\ge 2$) The moments on each internal face up to the degree $k-2$},\\
&\mbox{(for $k\ge 2$) The moments inside each element up to the degree $k-2$}.
\end{aligned}
\right.
\end{equation}

Then, apparently, one can follow, for a three-dimensional problem like \eqref{Pbmod}, the same path that we used for the two-dimensional case, {both for the {\it full Dirichlet} and the {\it full Neumann} case.} 

Given a decomposition $\Th$ of $\Omega$  into polyhedra $P$,  an integer $k\ge 1$,   and a smooth-enough function $g$ defined on $\Gamma$, we can define the finite dimensional subspaces   
$V_{h}^{g}$ and $V_{h}^{0}$ of $V_{h}(\Omega)$.

Still following slavishly the 2-dimensional path,  for each $P\in\Th$ and for each virtual element function $v$ we can define its $H^1_0(P)$-projection $\Pi_k^{\nabla}v$ as the unique solution (up to a constant that can be easily fixed) in $\Po_k$ of
\begin{equation}\label{defPinabla3}
\int_P\nabla (\Pi_k^{\nabla}v)\cdot\nabla p_k\dd P=\int_P\nabla v\cdot\nabla p_k\dd P \quad \forall p_k\in\Po_k(P).
\end{equation}
Ooops! When we attempt {\it to compute} the right-hand side of \eqref{defPinabla3} using the degrees of freedom \eqref{dof3d}, 
 we have
\begin{equation}\label{oops}
\int_P\nabla v\cdot\nabla p_k\dd P=\int_{\partial P} v \frac{\partial  p_k}{\partial n}\dd \sigma -\int_P v\,\Delta p_k\dd P.
\end{equation}
Now, the second term in the right-hand side of \eqref{oops} does not cause any trouble: for $p_k$ in $\Po_k$ we have that 
$\Delta p_k\in\Po_{k-2}$ and the term can be computed using the degrees of freedom \eqref{dof3d}. But for the first term
in the right-hand side of \eqref{oops} we would need to know the moments of $v$ on each face up to the order $k-1$ (the degree of $\frac{\partial  p_k}{\partial n}$), while in \eqref{dof3d} we have the moments only up to $k-2$. The way-out, as presented first in \cite{projectors}, is: 
\begin{equation}\label{trucco}
	\mbox{on each face $\ff$, replace }\int_{\ff} v \frac{\partial  p_k}{\partial n}\dd {\ff} \;\mbox{ with }\int_{\ff}\Pi^{\nabla,\ff}_k v \frac{\partial  p_k}{\partial n}\dd \ff,
\end{equation}
where $\Pi^{\nabla,\ff}_k v$  is the two-dimensional projection of $v$ onto $\Po_k(\ff)$, as defined in \eqref{defpinabla}-\eqref{fixc}, whose computation, in turn, on each face $\ff$  requires the moments of $v$ on $\ff$ only up to the order $k-2$. One might consider this as a typical use of the approach described in  subsection \ref{projectors}.

For the Neumann case, also the treatment of the right-hand side needs a more careful approach than in the 2D case. The integrals on $\Gamma_N$ might be computed as in \eqref{trucco}:
\begin{equation*}
\mbox{on each face in $ \Gamma_N$ replace} \quad \int_{\sigma\in \Gamma_N} g_N \,v \dd \sigma \quad \mbox{ with } \quad  
\int_{\sigma\in \Gamma_N}g_N \,\Pi^{\nabla,\sigma}_k v \dd \sigma .
\end{equation*}
Finally, the term $(f,v_h)$ can be treated as in \eqref{terminenoto}. For  more details we refer for instance to \cite{BBDMR-Cina}.

\begin{remark} As pointed out in subsection \ref{projectors}, we might prefer other projectors (instead of $\Pi^{\nabla}_k$) and use the \lq\lq{}projected $v$\rq\rq{} in place of $v$ for other quantities that cannot be computed directly from the degrees of freedom of $v$.
\qed
\end{remark}

Once $\Pi^{\nabla}_k$ has been defined, we can follow step-by-step (with obvious minor modifications) the path of the two-dimensional case (as suggested in \cite{projectors}) both for the {\it full Dirichlet} and the {\it full Neumann} case.

\begin{remark}\label{MFD-like}
Similarly to what was discussed and suggested in Subsection \ref{mix}, here too we could consider the possibility of using
other degrees of freedom, in particular on the ``Neumann boundary faces". This could be seen as: As usual the unknowns in
the computer are the degrees of freedom. To these degrees of freedom, in most elements (with few exceptions, as we shall see)
we could associate {\it a space of local functions}, some explicitly computable (e.g. polynomials), others  not explicitly computable but well defined, as solutions
of local PDE problems (that we will not solve, and use instead their projections on polynomial spaces). Other degrees of freedom might be used directly to compute projected polynomials, more in the spirit of  Mimetic Finite Differences. What we will carefully preserve will always be the {\it Patch test}, ensuring that the method will  be exact whenever the exact solution is a polynomial.  
\qed
\end{remark}

\section{2D VEM for  polygons with a curved edge}\label{2dcurvi}

Assume now that we have a  ``polygon"  $P$ with {\bf one} curved edge, that we call $\eta$, belonging to $\partial \Omega$. The case of polygons having two or more curved edges (always on $\partial \Omega$)
could be treated in a very similar manner. 

\subsection{The full Dirichlet case for 2D problems with curved boundary}\label{fullDir}

For the full Dirichlet problem, we can follow slavishly what has been done in \eqref{VkfP}-\eqref{VkfOm}, as well as in 
\eqref{defahpd}. Essentially, since the values of  the VEM spaces on the boundary are {\it assigned} (both for test and trial functions) the fact that the boundary is curved is just a minor nuisance, requiring the computation of integrals on the curve.
As we shall see, this will not be the case for the Neumann boundaries.

\subsection{\color{black}The full Neumann case for 2D problems with curved boundary}\label{farloc}

\begin{figure}[!htb]
\begin{center}
\includegraphics[height=0.30\textwidth]{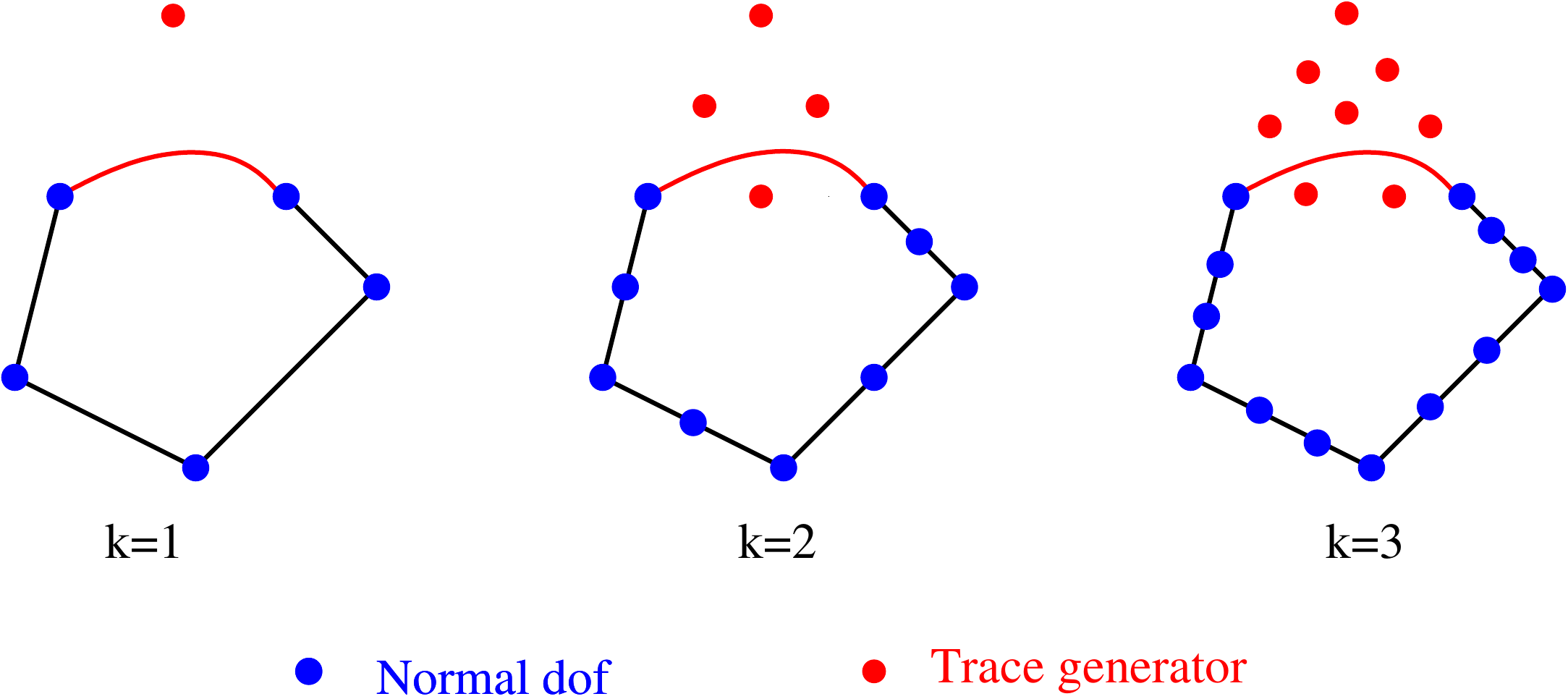} 
\end{center}
\caption{{{\color{blue} Dofs} }and {\color{red}Generating points} of $ \BB^{\eta}_k(\partial P)$ for  $k=1,2,3$}
\label{Farlocchi}
\end{figure}

We recall in this subsection the procedure introduced in \cite{Curvi2d}, which we refer to for more details.
In a natural way, on the curved edge $\eta$ we introduce
\begin{equation}\label{pokdieta} 
\Po_k(\eta):=\{v\in C^0(\eta)  \mbox{ s.t. } \exists\, p\in\Po_k(\R^2) \mbox{ with } v=p_{|\eta}\}.
\end{equation}
\begin{remark}\label{casino} Definition \eqref{pokdieta}, simple as it may seem, opens the door towards a more delicate discussion. Indeed, for a curved edge $\eta$ the dimension of the space $\Po_k(\eta)$ might be, in certain cases, far from obvious. In facts, the dimension could change between a minimum of $k+1$ (when the edge $\eta$ is straight), to a maximum of $k(k+1)/2$ when $\eta$ is {\it curved enough}
(here meaning that the only element in $\Po_k(\R^2)$ that vanishes identically on $\eta$ is the polynomial $\equiv 0$).  If we want to be allowed to treat systematically $\eta$ as {\it curved} we must then consider $k(k+1)/2$ generators (typically, the values at
$k(k+1)/2$ points as in Figure \ref{Farlocchi}), but the actual dimension might be as few as $k+1$. We shall come back to this problem (and related troubles)
pretty soon. \qed
\end{remark}

For the moment we can keep going on the same track as before, setting
\begin{equation*}
\begin{aligned}
\BB_k^{\eta}(\partial P)&:=
\{v\in C^0(\partial P) \mbox{ with } v_{|e} \in \Po_k(e)\, \forall \mbox{ straight edge } e \in \partial P, \\
&\hskip2cm \mbox{ and } v_{|\eta} \in \Po_k(\eta)\},
\end{aligned}
\end{equation*}
and then
\begin{equation}\label{localVkp}
V_k^{\eta}(P)=\{v\in H^1(P) \mbox{ such that } v_{|\partial P} \in \BB^{\eta}_k(\partial P)\\  \mbox{ and }\Delta v \in \Po_{k-2}(P)\}.
\end{equation}
So far so good. Now we observe that an element $v\in V^{\eta}_k(P)$ will be uniquely determined by:
\begin{align}
\label{dofs-cface-vertc}
&\bullet \mbox{the values of  }v \mbox{ at the vertices of  }P,\\
\label{dof-face-edgc}
&\bullet\mbox{(for }k\ge 2)\quad\int_e v\, p_{k-2}\dd e\quad\forall \mbox{ straight edge } e, \forall p_{k-2}\in \Po_{k-2}(e),\\
\label{dofs-cface-cedgc}
&\bullet\int_\eta v\, p_{k}^0\dd \eta\quad\forall p_{k}^0\in \Po_{k}(\R^2)\cap H^1_0(\eta),\\
\label{dofs-cface-insc}
&\bullet \mbox{(for }k\ge 2)\quad\int_P v\, p_{k-2}\dd P\quad\forall p_{k-2}\in \Po_{k-2}(P).
\end{align}

It is {\it crucial} to point out that, in \eqref{dofs-cface-cedgc}, we did not {\it forget} to say ``for $k\ge 2$". Indeed, as we already discussed in Remark \ref{casino},  for $k=1$,   {\it on a straight edge} $e$ we do have clearly that {\it all polynomials of degree $\le 1$ vanishing at the two endpoints of $e$ would be identcally zero on the whole edge}. But this is not anymore true for a {\it curved} edge, where a polynomial of degree $1$ in $\R^2$, vanishing at two distinct points would (indeed!) vanish identically on the straight line connecting the two points, but might very well be different from zero, for instance, on a circular arch passing through the same two points. And, in fact, {\it hic sunt leones} (that is:
{\it here come the difficulties}). Actually, when analysing and coding a numerical method, the basic step is often to define the finite dimensional space where we are looking for the solution, and the basis to be used for it (the {\it degrees of freedom}). The difficulty here  is in finding the dimension of such a space, or, in the terminology of Numerical Analysis, the correct  number of {degrees of freedom}. In particular, {\it to forbid   curved edges that are nearly flat} would be quite cumbersome in applications: the crucial point would be the definition of {\it nearly flat}: resonably clear for the human mind, but a nasty source of troubles in a computer code, taking into account that there the definition must be quantitative, and a minor change in it can produce a significant change in the numerical solution. And, sticking for simplicity to the case $k=1$ (but the difficulty pops up for every $k$), on a general curved edge $\eta$ the space  $\Po_k(\eta)$ has dimension 3, but when the edge $\eta$ is {\it almost} a straight line the number of parameters necessary to identify an element of  $\Po_k(\eta)$   {\it tends to} 2 (in a sense difficult to be made precise, but surely prospecting troubles).

\begin{figure}[!htb]
\begin{center}
\includegraphics[height=0.22\textwidth]{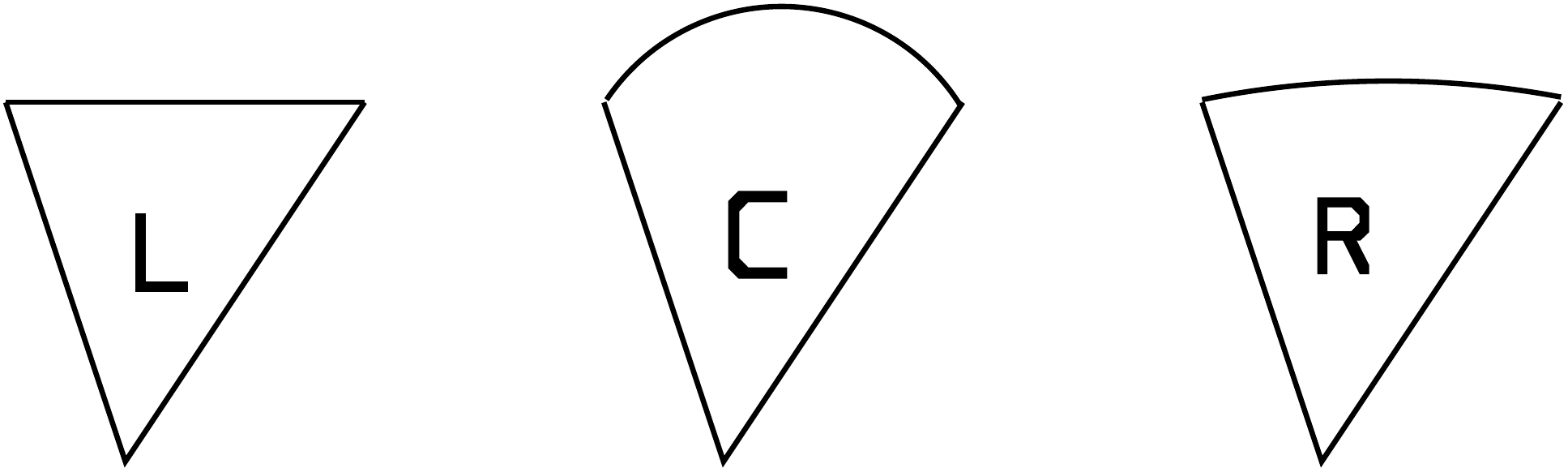} 
\end{center}
\caption{The dimension of $V_1(L)$ is 3. That of $V_1(C)$ is 4. And that of $V_1(R)$??}
\label{Different cases}
\end{figure}

 Referring to Figure \ref{Different cases}, for $k=1$ the dimension of $V_k(P)$ depends on the shape of the {\it curved edge}. In particular,
in {\it exact arithmetic} the dimension of $V_1(R)$ is clearly 4. But in the computer we should treat it carefully, since a simple minded treatment might lead to singular or nearly singular matrices depending on the number of digits that we are using.

The typical choice for VEMs (as done in \cite{Curvi2d})  is to treat every edge $\eta$ that {\it could  be curved} as {\it curved}, and to consider as unknowns 
{\it associated with the edge} the values (on that edge)  of  all polynomials of $\Po_k(\R^2)$ that vanish at the two endpoints of the edge, that therefore would sum up to $(k+1)(k+2)/2-2$ parameters. Typically, for a \lq\lq{}curved\rq\rq{} edge $\eta$ this is done considering  first the straight segment $Q$ connecting the two endpoints of $\eta$, then  considering an equilateral triangle $T_{\eta}$  having  $Q$ as one of its edges,  and finally  considering on $T_{\eta}$ the traditional degrees of freedom that one would have for $\Po_k$ on that triangle:
\begin{itemize}
\item the value at the vertex not belonging to $\eta$,
\item for $k\ge 2$ the values  at $k-1$ equally spaced points on each edge, and 
\item {for $k\ge 3$} the values at $(k-1)(k-2)/2$ internal nodes. 
	\end{itemize}
(These will be called {\bf generating points}).  
Hence, on an edge $\eta$ that has been declared as {\it curved} it will be simpler (in the computer code) to replace 
	\eqref{dofs-cface-cedgc} with the value at the generating points.

Therefore, even when the {\it so called ``curved edge"} is in fact straight, we will associate with it $(k+1)(k+2)/2-2$ parameters: the values at the {\it generating points} (or,  alternatively, the moments on $\eta$ against the traces (on $\eta$) of all $\Po_k$ polynomials vanishing at the two endpoints of $\eta$). Clearly this will be a set of {\it generators} (in the sense of \eqref{defGen}), but not a set of degrees of freedom, in the traditional sense.

All this, unfortunately, will require some additional paraphernalia.  Indeed, the obvious path, once we defined the local spaces
\eqref{localVkp}, is to collect them in the global space
\begin{equation}\label{globalVkp}
V_h\equiv V_h(\Omega):=\{v\in C^0(\overline\Omega)\mbox{ such that }v_{|P}\in V^{\eta}_k(P)\, \forall  P\in \Th\},
\end{equation}
where $\Th$ is our decomposition of $\Omega$ into \lq\lq{}polygonal\rq\rq{} elements, possibly with a curved edge. Then
 we would like to construct an approximate problem of the type \eqref{Pbmodh2d}. But in the computer code, the {\it unknowns} of the discretized problem cannot be the elements of $V_h$ (that are {\it functions}),
and {\bf must} be, instead, their {\it generators} (that are elements of some suitable $\R^N$). Remember that: given the generators, the function is uniquely determined, but not the other way around, as (for instance when the edge is {\it straight}) the same function on $\eta$ could be generated by several different  generators. Hence, once we have defined the global spaces in \eqref{globalVkp}, we must consider the space 
${\M}_{V_h}$
 ($\M$  is for {\it Midwife}!)  given by
\begin{equation*}
{\M}_{V_h} := \{\mbox{all the generators $\gv$ of the elements of $V_h$}\},
\end{equation*}
that, in general, might have a dimension bigger than that of $V_h$.

Hence we will have a linear mapping 
\begin{equation*}
\MM:\quad \M_{V_h}\;\longrightarrow V_h
\end{equation*} 
($\MM$ is for {\it Children}) that is, in other words:
\begin{equation*}
\MM\in{\mathcal L}(\M_{V_h},V_h). 
\end{equation*}

For every  given element $\gv\in\M_{V_h}$, $\MM$ selects a unique element $v=\MM(\gv)\in V_h$,  {\it generated by} $\gv$, that will be the function that we are interested in.

We observe that, given a generator  $\gv\in\M_{V_h}$, the polynomial  $\Pi^\nabla_k(\MM(\gv))$ can be computed as in \eqref{defpinabla}-\eqref{fixc} since the function $v=\MM(\gv)$ is defined all over the boundary of the element. Needless to say, other projectors from $V_k(P)$ to $\Po_k$  could be computed, if needed, for instance using the degrees of freedom as in  \eqref{proD}.

Then, with obvious notation, for every element $P$  we can define the bilinear form
\begin{equation}\label{defahpc}
a_h^P(\gu,\gv):=a^P(\Pi^\nabla_k({\MM}(\gu)),\Pi^\nabla_k({\MM}(\gv)))+\SS_P(\gu-{\mathcal D}(\Pi^{\nabla}_k({\MM}(\gu))),\gv-{\mathcal D}(\Pi^{\nabla}_k({\MM}(\gv)))),
\end{equation}
 where again $\SS_P$ is a bilinear form on $\R^N$ that scales like $a^P(\Pi^\nabla_k({\MM}(\gu)),\Pi^\nabla_k({\MM}(\gv)))$ and is positive on the kernel of ${\mathcal D}(\Pi^\nabla_k\MM)$. Note that, even though ${\mathcal D}$ (a right inverse of ${\MM}$  as in \eqref{defD}) is not uniquely defined (unless $ {\MM}$ is injective), ${\mathcal D}(\Pi^\nabla_k(\MM(\gv)))$ can be uniquely defined since $\Pi^\nabla_k(\MM(\gv))$
is a polynomial.

Then, as we did for polygons with straight edges, we can collect the definitions \eqref{defahpc} setting
\begin{equation*}
a_h(\gu,\gv):=\sum_{P\in \Th} a^P_h(\gu,\gv)\qquad \forall\, \gu,\gv\in \M_{V_h}.
\end{equation*}

The discretised problem reads now:
\begin{equation*}
\left\{
\begin{aligned}
&\mbox{Find $\gu\in \M_{V_h}$ such that}:\\
& a_h(\gu,\gv)=(f_h,\MM (\gv)) + <g_N,\MM (\gv)>_{\Gamma_N} \quad\forall \gv \in \M_{V_h} .
\end{aligned}
\right.
\end{equation*}
The treatment of the term $ (f_h,\MM (\gv))$ is discussed in \cite{Curvi2d}, which we refer to. Neumann boundary conditions were not dealt with in \cite{Curvi2d}, but the computation of the second term in the right-hand side does not pose additional difficulties.

\begin{remark}\label{curvinterni}
For problems where curved edges occur {\it inside} the domain $\Omega$, one has to be careful with the choice of {\it generators}. Indeed,  an edge that is internal to $\Omega$ will naturally belong to the boundary 
of two different elements. But the functions
$v\in V_h$ (as defined in \eqref{globalVkp}) must be {\it single valued} on the common edge, so that the set of generating points used must be {\it the same} for the two elements.  Moreover, as pointed out in \cite{Curvi2d}, 
it is better to stabilize the generators for the {\it curved common edge} only once: choosing, once and for all, one of the two elements having the edge in common, and then using them in the stabilizing term only when dealing with the chosen element. 
\qed
\end{remark} 

 The approach presented in the previous section does not extend easily to three-dimensional problems. In this subsection we introduce a couple of new ideas that could be used instead.

\subsection{Using a subset of degrees of freedom.}\label{interni}

Another way to set the values of VEM functions on curved edges of $\partial\Omega$ (in particular
on the part of $\partial\Omega$ where natural boundary conditions are prescribed) would be to
use the degrees of freedom on the straight edges plus the internal ones. To explain the procedure, let $P$ be a polygon, with just one curved edge $\eta$ belonging to $\partial\Omega$. We would like to define a space of functions like  \eqref{defBkf}-\eqref{Vkf}, that is, a space of continuous functions, polynomials of degree $k$ on each edge, and with Laplacian polynomial of degree $k-2$ in $P$. The problem is that we do not know how to define these functions on the curved edge. The information that we have (candidates to become degrees of freedom) are:
\begin{equation}\label{dof-subset}
\left\{
\begin{aligned}
&\mbox{the values of }v \mbox{ at the vertices of  }P,\\
&\mbox{(for }k\ge 2)\quad\int_e v\, p_{k-2}\dd e\quad\forall p_{k-2}\in \Po_{k-2}(e)\quad\forall\mbox{ {\bf straight} edge }e \subset  \partial P,\\
&\mbox{(for }k\ge 2)\quad\int_P v\, p_{k-2}\dd P\quad\forall p_{k-2}\in \Po_{k-2}(P).
\end{aligned}
\right.
\end{equation}
These conditions are not enough to individuate a VEM-function, but they are enough to individuate a polynomial of degree $k$.
Indeed, to make the simplest possible example, let us consider a triangle-like element,  i.e., with two straight edges and one, $\eta$, curved. For every integer $ k \ge 1$ we see that using the values of $v$ on the two straight
edges (amounting to  $ 2k + 1 $ dofs), and the internal moments up to the degree $k-2$ (amounting to $k(k-1)/2$ dofs) we can identify uniquely a polynomial $p_k$ in $\Po_k$. Indeed,
\begin{equation*} 
2k+1+k(k-1)/2=(k+1)(k+2)/2
 ~\equiv dim(\Po_k).
\end{equation*}
With a given VEM-function $v$ we can associate a polynomial $p_k^* \in \Po_k$, computed  through the conditions
\begin{equation}\label{dof-subset-pk}
\left\{
\begin{aligned}
&~p_k^*= v \mbox{ at the vertices of  }P, \quad \mbox{ and for } k\ge 2,\\
&\int_e p_k^*\, p_{k-2}\dd e=\int_e v\, p_{k-2}\dd e ~\forall p_{k-2}\in \Po_{k-2}(e)\quad\forall\mbox{ {\bf straight} edge }e \subset  \partial P,\\
&\int_P p_k^*\, p_{k-2}\dd P=\int_P v\, p_{k-2}\dd P\quad\forall p_{k-2}\in \Po_{k-2}(P).
\end{aligned}
\right.
\end{equation}
It is trivial to check that (always in the case of a triangle-like element) conditions \eqref{dof-subset-pk} are a set of unisolvent degrees of freedom for $\Po_k$, so that $p_k^*$ is unique.
Then we can take 
\begin{equation}\label{dof-subset-eta}
v_{|\eta}=p^*_{k|\eta}.
\end{equation}
Now conditions \eqref{dof-subset} plus  \eqref{dof-subset-eta} determine $v$ uniquely, and we
 point out that, by construction, global continuity is guaranteed. We observe that, in this particular case of triangle-like elements, the function $v$ is actually a polynomial, coinciding with $p_k^*$ on the whole element.

This result could suggest to have a decomposition made of triangle-like around the boundary. The approach would be simple and clean, but not really in the spirit of VEM, whose strong point is to allow polygons of arbitrary shape.

Let us then consider a general polygon, with a number of straight edges higher than 2. Then we would have more information than necessary, i.e., the number of conditions \eqref{dof-subset} would be bigger than the dimension of $\Po_k$. In this case
we should use a least-square solution, keeping of course fixed the values at the two endpoints of $\eta$ to guarantee the global continuity. More precisely, if $N$ is the number of conditions \eqref{dof-subset}, ordered in such a way that the values at the two endpoints of $\eta$ are the last two, we solve the following problem:
\begin{equation}\label{least-square}
\left\{
\begin{aligned}
&\mbox{find } p_k^* \in \Po_k(P) \mbox{ such that } p_k^*=v \mbox{ at the two endpoints of } \eta, \mbox{ and}\\
& \sum_{i=1}^{N-2} (dof_i(p_k^*)-dof_i(v))^2 = \mbox{ minimum}.
\end{aligned}
\right.
\end{equation}
Then, as before, we  define 
\begin{equation}\label{dof-subset-eta-gen}
v_{|\eta}=p^*_{k|\eta}.
\end{equation}
Once the function $v$ is individuated by \eqref{dof-subset} and \eqref{dof-subset-eta-gen}, we can define  $\Pi^{\nabla}_k v$ as in \eqref{defpinabla} (or $\Pi^0_{k-1}\nabla v$  as in \eqref{projgrad}) and write the discrete bilinear form as in \eqref{defahpd} (or as in \eqref{variable}, respectively).

\begin{remark}
This strategy cannot be used, as such, if the curved
edge is inside $\Omega$ and common to two or more elements. In this case one might think of a
{\it master and slave} approach (only one element sets the degrees of freedom to be used on
the curved edge) or take a suitable combination of the effects of the elements sharing the
curved edge. Once the functions are defined on the curved edges, everything goes along the same lines used  for  polygons with straight edges, and the discrete problem can be written exactly as in \eqref{Pbmodh2d-Neumann}.
\qed
\end{remark}

\begin{remark} \label{nofunz}A similar approach, although conceptually very different, would consist in using, associated with the curved edge, suitable
{\it degrees of freedom} to be used in the construction of the various projectors, {\it without defining a functional space within the element}. In other words, in the element with a curved edge we would consider only {\it degrees of freedom} without introducing a functional space. Hence, associated with the element with the curved edge we will still have a local stiffness matrix (as we would have
when using {\it Mimetic Finite Differences}). This would correspond, somehow, to use Finite Elements (or Virtual Elements) in the other elements, and Mimetic Finite Differences in the elements having  a curved edge with natural boundary conditions (in the spirit, somehow, of Subsection \ref{mix}). Although this might look as a {\it monster} from the theoretical point of view, the computer code will not suffer (in particular if you already have a VEM code and an MFD one) as, in any case, it will deal only with degrees of freedom and never with functions. The only drawbacks would appear only in the proof of error estimates, although, as a combination of two reliable methods, it is reasonable to expect the usual level of  accuracy.
The above procedure would be in this case the following: we use again a least-square approach to compute a polynomial $p^*_k\in \Po_k$,  this time without any need for fixing its values at the two endpoints. Then, to individuate a VEM-function $v$ we add to conditions \eqref{dof-subset} the moments of $p^*_k$ of order $k-1$ as degrees of freedom on $\eta$:
\begin{equation}\label{dof-eta-MFD}
\int_{\eta} v\, \nabla p_{k}\cdot {\bf n} \dd \eta=\int_{\eta} p^*_k\, \nabla p_{k}\cdot {\bf n}\dd \eta\quad\forall p_{k}\in \Po_{k}(\eta).
\end{equation}
Conditions \eqref{dof-subset} and \eqref{dof-eta-MFD} allow us to compute the projection  $\Pi^{\nabla}_k v$ as in \eqref{defpinabla} (or $\Pi^0_{k-1}\nabla v$  as in \eqref{projgrad}), and we can write the discrete bilinear form as in \eqref{defahpd} (or as in \eqref{variable}, respectively). 
\qed
\end{remark} 

\subsection{Coupling FEM and VEM: the superimposed polygon}\label{theribbon}
Another possible approach that might be considered (always on the part of the boundary where Neumann boundary conditions are imposed) is the use of a {\it superimposed polygon}
{ $\Pit$}  including  {\color{blue}$\Omega$}, such that the two  boundaries ($\partial\Pit$ and $\partial \Omega$) are never closer to each other more than $C h$ for some fixed $C>0$. 
\begin{figure}[htp]
\begin{center}
\includegraphics[height=0.62\textwidth]{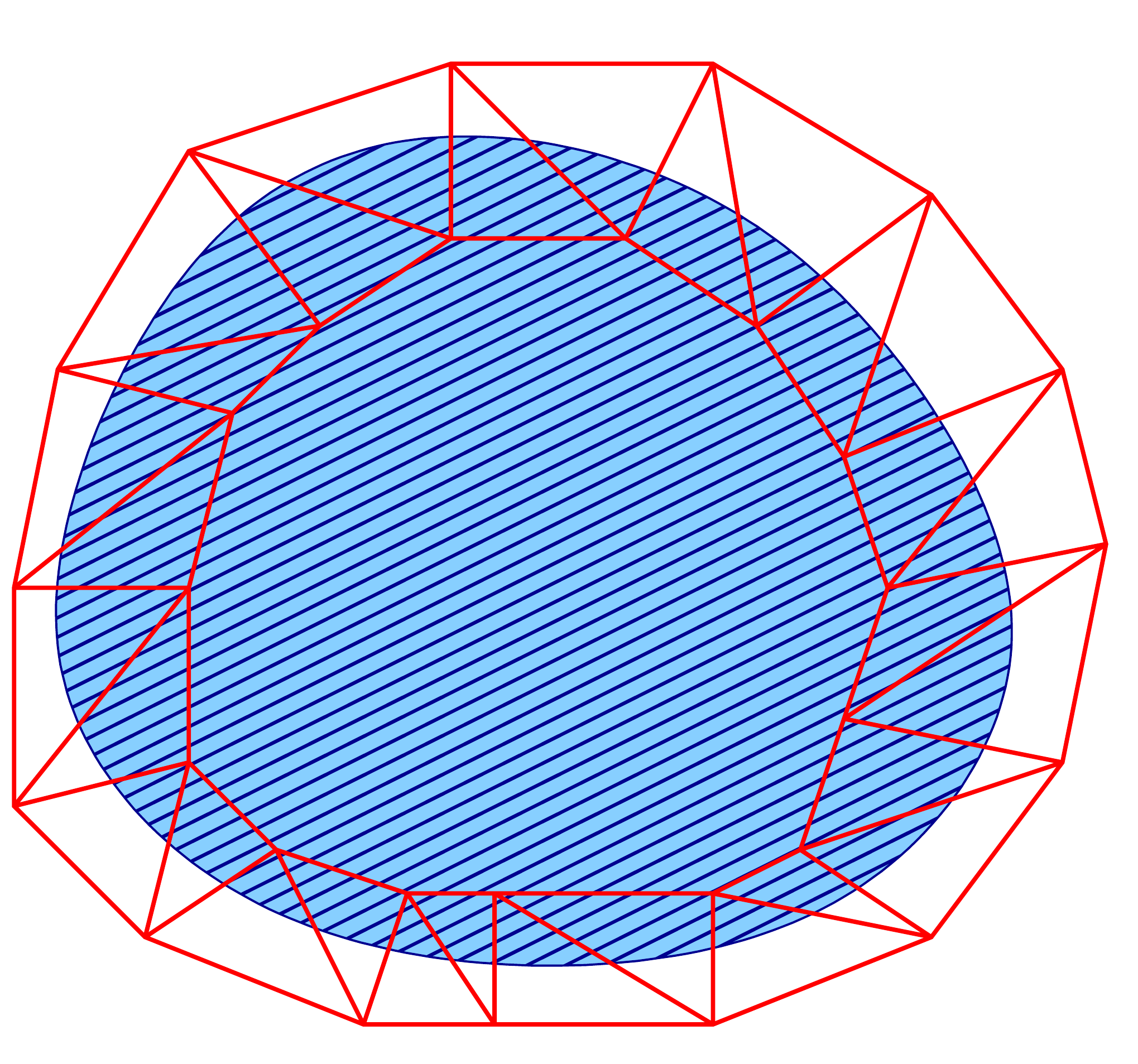} 
\end{center}
\caption{{\color{blue} $\Omega$} and {\color{red}the ribbon}}
\label{Corona}
\end{figure}
 Let then $\TTT(\Pit)$ be a decomposition  of $\Pit$ in polygons  (for simplicity, in convex polygons), that naturally produces a decomposition  of  $\Omega$  consisting of {\it normal polygons} (corresponding to the elements of $\TTT(\Pit)$ that are all contained in $\Omega$) and {\it polygons with a curved edge} (corresponding to the restrictions to $\Omega$ of the polygons  in $\TTT(\Pit)$  that contain parts of $\Omega$).  For simplicity we assume that there are no polygons in $\TTT(\Pit)$  that do not contain at least a part of $\Omega$. 

A simple way of realizing this is to construct a ribbon of quadrilaterals, as in Fig.\ref{Corona}, around $\partial\Omega$, that naturally defines { two polygons: one {\it containing} $\Omega$ (that will be our $\Pit$), and one (say, $\widetilde Q$), {\it contained in} $\Omega$. The polygon $\widetilde Q$ can then be decomposed as we need, while the quads in the ribbon will be decomposed in triangles. Clearly, the union of the decomposition of $\widetilde Q$ and the strip of triangles forms a decomposition of $\Pit$.}

Once $\Pit$ and its decomposition have been constructed, we can define the discrete space
 ${V}_h(\Pit)$ as follows. On polygons belonging to $\widetilde Q$ we will use virtual elements of degree $k$, as we did in \eqref{Vkf}, while on the triangles in the ribbon we simply take usual finite elements of degree $k$.
  Notice that global continuity is ensured, since the degrees of freedom at the interelements are the same for Virtual and Finite elements. The unknowns of our problem will then be the elements of ${V}_h(\Pit)$. The next step is to construct a suitable projection operator $\Pi^{\nabla}_k$ in each element of the decomposition. For elements in $\widetilde Q$ which, according to the construction, will be just  polygons as those that we considered before, we proceed as in the previous Section \ref{2dflat}. The discrete bilinear form  will be defined as in \eqref{defahpd} (or in \eqref{variable})  in the elements contained in $\widetilde Q$, after computing the $\Pi^{\nabla}_k$-operator as in \eqref{defpinabla} (or $\Pi^0_{k-1}\nabla v$  as in \eqref{projgrad}).

For the triangles in the ribbon we do not need any projection, since we already have polynomials. Hence, on the {\it true} element
\begin{equation*}
 \Omega_{T}:=\Omega\cap T 
\end{equation*}
we can simply take $\Pi^{\nabla, \Omega_{T}}_k v={p_k}_{| \Omega_{T}}$.
Then, the discrete bilinear form will be defined simply as \eqref{defauv} in each  $\Omega_T$.

\section{3d VEM for  polyhedra with a curved face}\label{3dcurvi}

The extension of the approach described in Subsection \ref{farloc} to three-dimensional problems looks particularly hard. On the one hand, the number of generating points grows in a significant way with the degree $k$. On the other hand, most important, ensuring global continuity looks particularly complicated to realize. Instead, the new approaches indicated in Subsections \ref{interni} and \ref{theribbon} look more promising. Clearly, extensive numerical tests will need to be performed to validate them in terms of feasibility, efficiency, and accuracy. We briefly sketch here the extension of the two approaches. Here too, like we said in subsection \ref{fullDir}, the full Dirichlet case presents no difficulties. We then concentrate on the treatment of the full Neumann case, or of the part of the boundary where Neumann conditions are imposed.

\subsection{Using a subset of degrees of freedom.}\label{interni3d}

The extension of the procedure highlighted in subsection \ref{interni} requires some care. Actually, on a polyhedron the number of degrees of freedom not positioned on $\partial \Omega$ is in general higher than the dimension of $\Po_k$. Thus, a least-square solution would be necessary, but this is in contrast with the continuity requirement at the interfaces. In 2d the problem was circumvented easily just by fixing the values at the two endpoints of the curved edge. In 3d we should fix the value on the whole boundary of the curved face, which is what we are trying to define through the construction of the polynomial. Hence, in 3d we have only two possibilities, and not three like in 2d.

One possibility is to have a decomposition of tetrahedra-like around the boundary, that is, polyhedra with 3 flat faces and one, say $\sigma$, curved face on $\partial \Omega$. We see that the values of a VEM-function $v$ on $\sigma$ could be obtained through its dofs on the three flat faces plus the internal moments up to the order $k-3$ {\bf only} (and not $k-2$). More precisely,
 for every integer $ k \ge 1$, using the values of $v$ on the three straight
edges (amounting to  $ 3k + 1 $ dofs), the moments up to the degree $k-2$ on the three flat faces (amounting to $3 k(k-1)/2$ dofs), and the internal moments up to the degree $k-3$ (amounting to $k(k-1)(k-2)/6$ dofs), we can identify uniquely a polynomial $p^*_k\in \Po_k$. Indeed we have:
\begin{equation*}
(3k + 1)+\frac{3 k(k-1)}{2}+\frac{k(k-1)(k-2)}{6}= \frac{k^3+6k^2+11k+6}{6}\equiv  dim (\Po_k).
\end{equation*}
We can then compute a polynomial $p_k^* \in \Po_k(P)$ through the conditions
\begin{equation}\label{dof3d-pk}
\left\{
\begin{aligned}
&p_k^*=v \mbox{ at the vertices of }P,\\
&\mbox{for $k\ge 2$} ~ \displaystyle{\int_e p_k^*\, p_{k-2}\dd e\!\!=\!\!\int_e v\, p_{k-2}\dd e} ~ \forall p_{k-2}\in \Po_{k-2}(e) ~ \forall  \mbox{ {\bf straight} edge } e,\\
&\mbox{for $k\ge 2$} ~ \displaystyle{\int_{\ff} p_k^*\, p_{k-2}\dd {\ff}\!\!=\!\!\int_{\ff} v\, p_{k-2}\dd {\ff}} ~ \forall p_{k-2}\in \Po_{k-2}({\ff}) ~ \forall \mbox{ {\bf flat} face }\ff,\\
&\mbox{for $k\ge 3$} ~ \displaystyle{\int_{P} p_k^*\, p_{k-3}\dd {P}\!\!=\!\!\int_{P} v\, p_{k-3}\dd {P}} ~ \forall p_{k-3}\in \Po_{k-3}(P).
\end{aligned}
\right.
\end{equation}
(It is not difficult to check that the dofs \eqref{dof3d-pk} are unisolvent for $\Po_k$).
Once the polynomial $p^*_k$  has been computed, we take $v_{|\sigma}=p^*_{k|\sigma}$, and proceed as we did for polyhedra with flat faces. We point out that, by construction, the global continuity is guaranteed. We also observe that, contrary to what happens in 2D, the function $v$ does not coincide with the polynomial, since for computing $p^*_k$ we did not use all the internal moments of $v$, but just the moments up to order $k-3$.
To summarize, a function $v$ is completely determined by the following conditions:
\begin{equation}\label{dof3d-curved}
\left\{
\begin{aligned}
&\mbox{The values at the vertices of } P,\\
&\mbox{for $k\ge 2$  the moments}  \displaystyle{\int_e v\, p_{k-2}\dd e}~\forall p_{k-2}\in \Po_{k-2}(e)~\forall  \mbox{ {\bf straight} edge } e,\\
&\mbox{for $k\ge 2$  the moments } \displaystyle{\int_{\ff} v\, p_{k-2}\dd {\ff}}\quad\forall p_{k-2}\in \Po_{k-2}({\ff})\quad\forall \mbox{ {\bf flat} face }\ff ,\\
&\mbox{for $k\ge 2$ the moments } \displaystyle{\int_{P} v\, p_{k-2}\dd {P}}\quad\forall p_{k-2}\in \Po_{k-2}(P),
\end{aligned}
\right.
\end{equation}
\begin{equation}\label{dof-sigma}
v_{|\sigma}=p^*_{k|\sigma}.
\end{equation}
Conditions \eqref{dof3d-curved}-\eqref{dof-sigma} determine $v$ uniquely, allowing to compute $\Pi^{\nabla,\ff}_k v$ on each flat face $\ff$, and then $\Pi^{\nabla,P}_k v$ as in Section \ref{3dpiatti} (see  \eqref{defPinabla3}--\eqref{trucco}), so that the discrete bilinear form \eqref{defahpd} is computable. 

On a generic polyhedron, with more than three flat faces, the number of conditions \eqref{dof3d-pk} is bigger than the dimension of $\Po_k$. In such a case we might use a least-square procedure to compute a polynomial $p^*_k$, and then use its moments on $\sigma$ as degrees of freedom on the curved face. However, in order to guarantee the global continuity, before computing $p^*_k$ we have to take care of the flat faces.

\noindent
$\bullet 1^{st}$ step: on each flat face $\ff$ with one curved edge $\eta$ we use least-squares to compute a polynomial  $p^{\ff}_k\in \Po_k(\ff)$ as in \eqref{least-square}, and set $v_{|\eta}=p^{\ff}_{k|\eta}$, thus ensuring continuity of $v$ on $\partial \sigma$.

\noindent
$\bullet 2^{nd}$ step: we compute a polynomial $p^{*}_k\in \Po_k(P)$ as the least-square solution of \eqref{dof3d-pk}

\noindent
$\bullet 3^{rd}$ step: we set 
$$\int_{\sigma} v\,\nabla p_{k}\cdot {\bf n} \dd\sigma =\int_{\sigma} p^*_k\,\nabla p_{k}\cdot {\bf n} \dd\sigma ,\quad \forall p_{k}\in \Po_{k}(\sigma).$$
This procedure, together with conditions \eqref{dof3d-curved}, determines uniquely $v$, and allows to compute $\Pi^{\nabla,\ff}_k v$ on each flat face $\ff$, and then $\Pi^{\nabla,P}_k v$ as in Section \ref{3dpiatti} (see  \eqref{defPinabla3}--\eqref{trucco}), so that the discrete bilinear form \eqref{defahpd} is computable. 
Contrary to the case of tetrahedra-like elements, here we do not assign $v$ on the curved face. Instead we impose degrees of freedom, more in the spirit of Mimetic Finite Differences (see always Remark \ref{nofunz}).

\subsection{Coupling FEM and VEM: the superimposed polyhedron}\label{theribbon3}

Let us sketch how to extend  to the 3-D case  what we did in  Subsection \ref{theribbon}, extending to the {\it skin}  what we did in the two-dimensional case for  {\it the ribbon}.
For this, we assume that we are given a polyhedron  $\Pit$ that contains strictly our three-dimensional domain $\Omega$, and that $\Pit$ is (suitably) decomposed in polyhedra. More precisely, we will have a {\it skin} (that we manage to have  made of tetrahedra), while the {\it interior} of $\Pit$ will be decomposed into polyhedra as needed. 
In particular here too we will assume that all {\it the elements of the decomposition of $\Pit$ that are {\it internal to} $\Pit$} are also {\it internal to $\Omega$}, as we had in the 2-dimensional case. 
Then we can proceed, {\it mutatis mutandis}, as we did in the two-dimensional case. Once $\Pit$ and its decomposition have been constructed, we can define the discrete space
 ${V}_h(\Pit)$: we will use virtual elements of degree $k$, as we did in Section \ref{3dpiatti}, in the internal polyhedra,  while on the tetrahedra in the skin we simply take finite elements of degree $k$. The elements of ${V}_h(\Pit)$ will be our unknowns. The discrete bilinear form  will be defined as in \eqref{defahpd} (or in \eqref{variable})  in the elements contained in $\Omega$, after computing the $\Pi^{\nabla}_k$-operator as in Section \ref{3dpiatti}. In the tetrahedra, having already polynomials, we do not need any projection. Hence, like in the 2d case, in the elements  $\Omega_{T}:=\Omega\cap T$ with one curved face we take  $\Pi^{\nabla, \Omega_{T}}_k v_h= {p_k}_{| \Omega_{T}}$. Then, in each $\Omega_{T}$ the discrete bilinear form will simply be \eqref{defauv} as for Finite Elements.


\begin{thebibliography}{10}

\bibitem{projectors}
B.~Ahmad, A.~Alsaedi, F.~Brezzi, L.~D. Marini, and A.~Russo, \emph{Equivalent
  projectors for virtual element methods}, Comput. Math. Appl. \textbf{66}
  (2013), no.~3, 376--391.

\bibitem{Ovall}
A.~Anand, J.~S. Ovall, S.~E. Reynolds, and S.~Wei\ss~er, \emph{Trefftz finite
  elements on curvilinear polygons}, SIAM J. Sci. Comput. \textbf{42} (2020),
  no.~2, A1289--A1316.

\bibitem{ABCM}
D.N. Arnold, F.~Brezzi, B.~Cockburn, and L.D. Marini, \emph{Unified analysis of
  discontinuous {G}alerkin methods for elliptic problems}, SIAM J. Numer. Anal.
  \textbf{39} (2001), no.~5, 1749--1779.

\bibitem{Babuska71}
I.~Babu\v{s}ka, \emph{The rate of convergence for the finite element method},
  SIAM J. Numer. Anal. \textbf{8} (1971), 304--315.

\bibitem{Curvi2d}
L.~Beir\~{a}o~da Veiga, F.~Brezzi, L.~D. Marini, and A.~Russo, \emph{Polynomial
  preserving virtual elements with curved edges}, Math. Models Methods Appl.
  Sci. \textbf{30} (2020), no.~8, 1555--1590.

\bibitem{BLMbook}
L.~Beir\~{a}o~da Veiga, K.~Lipnikov, and G.~Manzini, \emph{The mimetic finite
  difference method for elliptic problems}, MS\&A. Modeling, Simulation and
  Applications, vol.~11, Springer, Cham, 2014.

\bibitem{rozzi}
L.~Beir\~{a}o~da Veiga, A.~Russo, and G.~Vacca, \emph{The virtual element
  method with curved edges}, ESAIM Math. Model. Numer. Anal. \textbf{53}
  (2019), no.~2, 375--404.

\bibitem{volley}
L.~Beir{\~a}o~{da Veiga}, F.~Brezzi, A.~Cangiani, G.~Manzini, L.D. Marini, and
  A.~Russo, \emph{Basic principles of virtual element methods}, Math. Models
  Methods Appl. Sci. \textbf{23} (2013), no.~1, 199--214.

\bibitem{BBDMR-Cina}
L.~Beir{\~a}o~{da Veiga}, F.~Brezzi, F.~Dassi, L.D. Marini, and A.~Russo,
  \emph{Serendipity virtual elements for general elliptic equations in three
  dimensions}, Chinese Annals of Mathematics Series B \textbf{39} (2018),
  no.~2, 315--334.

\bibitem{hitchhikers}
L.~Beir{\~a}o~{da Veiga}, F.~Brezzi, L.~D. Marini, and A.~Russo, \emph{The
  hitchhiker's guide to the virtual element method}, Math. Models Methods Appl.
  Sci. \textbf{24} (2014), no.~8, 1541--1573.

\bibitem{SERE-nod}
L.~Beir{\~a}o~{da Veiga}, F.~Brezzi, L.D. Marini, and A.~Russo,
  \emph{Serendipity nodal {VEM} spaces}, Comp. Fluids \textbf{141} (2016),
  2--12.

\bibitem{variable-primal}
L.~Beir{\~a}o~{da Veiga}, F.~Brezzi, L.D. Marini, and A.~Russo,
\emph{Virtual element methods for general second order elliptic
  problems on polygonal meshes}, Math. Models Methods Appl. Sci. \textbf{26}
  (2016), no.~4, 729--750.

\bibitem{VEM-curv-Bertol}
S.~Bertoluzza, M.~Pennacchio, and D.~Prada, \emph{High order vems on curved
  domains}, Rend. Lincei Mat. Appl. \textbf{30} (2019), no.~2, 391--412.

\bibitem{Brezzi:Marini:plates}
F.~Brezzi and L.D. Marini, \emph{Virtual element methods for plate bending
  problems}, Comput. Methods Appl. Mech. Engrg. \textbf{253} (2013), 455--462.

\bibitem{Chinosi:Marini:plates}
C.~Chinosi and L.D. Marini, \emph{Virtual element method for fourth orded
  problems: $l^2-$estimates}, Comp. Math. Appl. \textbf{76} (2016), 1959--1967.

\bibitem{cockburn-IMU}
B.~Cockburn, \emph{The hybridizable discontinuous {G}alerkin methods},
  Proceedings of the {I}nternational {C}ongress of {M}athematicians. {V}olume
  {IV}, Hindustan Book Agency, New Delhi, 2010, pp.~2749--2775.

\bibitem{TJR-Iso}
J.~A. Cottrell, T.~J.~R. Hughes, and Y.~Bazilevs, \emph{Isogeometric analysis},
  John Wiley \& Sons, Ltd., Chichester, 2009, Toward integration of CAD and
  FEA.

\bibitem{DiPietro-Ern}
D.~Di~Pietro and A.~Ern, \emph{Mathematical aspects of discontinuous {G}alerkin
  methods}, Math\'ematiques \& Applications (Berlin) [Mathematics \&
  Applications], vol.~69, Springer, Heidelberg, 2012.

\bibitem{Dole-Feista}
V.~Dolej{\v s}{\' i} and M.~Feistauer, \emph{Discontinuous {G}alerkin method.
  analysis and applications to compressible flow}, Springer Series in
  Computational Mathematics, vol.~48, Springer, Cham, 2015.

\bibitem{Gong-Wu-Xu}
S.~Gong, S.~Wu, and J.~Xu, \emph{New hybridized mixed methods for linear
  elasticity and optimal multilevel solvers}, Numer. Math. \textbf{141} (2019),
  no.~2, 569--604.

\bibitem{Lions:Magenes:1968}
J.-L. Lions and E.~Magenes, \emph{Probl\`emes aux limites non homog\`enes et
  applications. {V}ol. 1}, Travaux et Recherches Math\'ematiques, No. 17,
  Dunod, Paris, 1968.

\bibitem{Mu-Wang-Ye}
L.~Mu, J.~Wang, Y.~Wang, and X.~Ye, \emph{A computational study of the weak
  {G}alerkin method for second-order elliptic equations}, Numer. Algorithms
  \textbf{63} (2013), no.~4, 753--777.

\end{thebibliography}
\end{document}